\documentclass{amsart}
\usepackage{amsmath,amssymb,amsfonts}
\usepackage[mathscr]{eucal}

\usepackage{enumerate}

\input xy
\xyoption{all}

\theoremstyle{plain}
\newtheorem{theorem}{Theorem}[section]

\newtheorem{prop}[theorem]{Proposition}

\theoremstyle{definition}
\newtheorem{definition}[theorem]{Definition}
\newtheorem{example}[theorem]{Example}

\numberwithin{equation}{section}

\newcommand{\opname}[1]{\operatorname{\mathsf{#1}}}

\newcommand{\Mod}{\opname{Mod}\nolimits}

\newcommand{\per}{\opname{per}\nolimits}

\newcommand{\Ho}{\opname{Ho}\nolimits}
\newcommand{\id}{\text{id}}

\newcommand{\rep}{\opname{rep}\nolimits}

\newcommand{\colim}{\opname{colim}\nolimits}

\renewcommand{\ker}{\opname{ker}\nolimits}
\newcommand{\Hom}{\opname{Hom}}

\title{Cluster categories for topologists}
\author[J.\, E.\, Bergner]{Julia E.\ Bergner}
\address{Department of Mathematics, University of California, Riverside}
\email{bergnerj@member.ams.org}

\author[M.\, Robertson]{Marcy Robertson}
\address{Department of Mathematics, University of Western Ontario, Canada}
\email{mrober97@uwo.ca}

\keywords{triangulated categories, cluster categories, differential graded categories, stable model categories}

\begin{document}
\begin{abstract}
We consider triangulated orbit categories, with the motivating example of cluster categories, in their usual context of algebraic triangulated categories, then present them from another perspective in the framework of topological triangulated categories.
\end{abstract}
\maketitle
\tableofcontents

\section{Introduction}

Cluster algebras were introduced and studied by Berenstein, Fomin, and Zelevinsky \cite{fz1}, \cite{fz2}, \cite{fz3}, \cite{bfz}. It was the discovery of Marsh, Reineke, and Zelevinsky that they are closely connected to quiver representations \cite{mrz}.  This connection is reminiscent of one between quantum groups and quiver representations discovered by Ringel \cite{ringel} and investigated by many others.  The link between cluster algebras and quiver representations becomes especially beautiful if, instead of categories of quiver representations, one considers certain triangulated categories deduced from them. These triangulated categories are called \emph{cluster categories}.

Cluster categories were introduced by Buan, Marsh, Reineke, Reiten, and Todorov in \cite{bmrrt} and, for Dynkin quivers of type $A_n$, in the paper of Caldero, Chapoton, and Schiffler \cite{ccs}. If $k$ is a field and $Q$ a quiver without oriented cycles, the associated cluster category $\mathcal C_Q$ is the ``largest" $2$-Calabi-Yau category under the derived category of representations of $Q$ over $k$. This category fully determines the combinatorics of the cluster algebra associated with $Q$ and, simultaneously, carries considerably more information which was used to prove significant new results on cluster algebras.

The goal of this paper is to introduce topological triangulated orbit categories, and in particular, the motivating example of topological cluster categories.  In doing so, we hope to explain the fundamental ideas of triangulated orbit categories to readers from a more homotopy-theoretic, rather than algebraic, background.  Our goal, as topologists, is to work the theory of cluster categories backward, by understanding triangulated categories that have similar properties but which arise from purely topological origins. In particular, in future work we aim to provide sufficient conditions on a stable model category (or more general cofibration category) $\mathcal C$, equipped with a self equivalence $F \colon \mathcal C \rightarrow \mathcal C$, so that the orbit category $\mathcal C /F$ admits a triangulated structure.

We begin by presenting the definition of triangulated orbit categories in Section 2.  In Section 3, we elaborate on the notion of algebraic triangulated category and discuss the enhanced version of orbit categories in differential graded categories.  We conclude that section with a brief introduction to cluster categories, the primary example of interest.  In Section 4, we introduce topological triangulated categories and give definitions of topological triangulated orbit categories and the corresponding example of cluster categories.

\section{Triangulated orbit categories}

\begin{definition}
Let $\mathcal T$ be an additive category and $F \colon \mathcal T \rightarrow \mathcal T$ a self-equivalence of $\mathcal T$.  The \emph{orbit category} of $\mathcal T$ by $F$ is the category $\mathcal T/F$ with objects those of $\mathcal T$ and morphisms defined by
\[ \Hom_{\mathcal T/F}(X,Y) = \bigoplus_{n \in \mathbb Z} \Hom_\mathcal T(X, F^nY). \]
\end{definition}

The composite of a morphism  $f \colon X \rightarrow F^{n}Y$ with a morphism $g \colon Y\rightarrow F^{p}Z$ is given by $(F^{n}g)\circ f$.

Although the orbit category in fact has many more morphisms than the original category $\mathcal T$, we regard is as a kind of quotient; in particular, it comes equipped with a ``projection" functor $\pi \colon \mathcal T \rightarrow \mathcal T/F$, together with an equivalence of functors $\pi \circ F \rightarrow \pi$ which is universal with respect to all such functors.

If we merely require $\mathcal T$ to be an additive category, it is not hard to see that the orbit category $\mathcal T/F$ is again an additive category, and the projection $\pi \colon \mathcal T \rightarrow \mathcal T/F$ is an additive functor.  However, we are most interested in the case where $\mathcal T$ is in fact a triangulated category.  The question of whether $\mathcal T/F$ still has a natural triangulated structure is much more difficult.

Most basically, we would like to complete any morphism $X \rightarrow Y$ in $\mathcal T/F$ to a distinguished triangle.  If it comes from a morphism $X \rightarrow Y$ in $\mathcal T$, then there is no problem.  However, in general, it is of the form
\[ X \rightarrow \bigoplus_{i=1}^N F^{n_i}Y \]
in terms of maps in $Y$, and in this case is not clear how to complete such a morphism to a triangle in the orbit category.

In \cite{keller}, Keller gives conditions under which the orbit category associated to some algebraic triangulated categories still possess a natural triangulated structure.  He constructs a triangulated category into which the orbit category embeds, called the \emph{triangulated hull}, then shows under which hypotheses this triangulated hull is in fact equivalent to the orbit category.  While his conditions are fairly restrictive, he shows that they hold in several important applications.  Most significantly, they hold for the construction of the cluster category.

\section{Algebraic triangulated categories}

A triangulated category is \emph{algebraic} if it admits a differential graded model, sometimes referred to as an \emph{enhanced} algebraic triangulated category.

We deviate from algebraists' standard conventions in two minor points.  First, in line with grading conventions in topology, we grade complexes homologically (as opposed to cohomologically), so that differentials \emph{decrease} the degree by~1. Second, we use covariant (as opposed to contravariant) representable functors; the resulting dg categories we obtain are hence the opposite of those obtained dually.

Let $k$ be a field. A \emph{differential graded category}, or simply \emph{dg category}, is a category $\mathcal C$ enriched in chain complexes of $k$-modules. In other words, a dg category consists of a class of objects together with a complex $\Hom_{\mathcal C}(X,Y)$ of morphisms for every pair of objects $X,Y$ in $\mathcal C$. Composition is given by the tensor product of chain complexes, i.e.,
\[ \circ \colon \Hom_{\mathcal C} (Y,Z) \otimes_k \Hom_{\mathcal C}(X,Y) \to \Hom_{\mathcal C}(X,Z) \]
for all $X,Y,Z$ in $\mathcal C$ which is associative and admits two-sided units $1_X \in \Hom_{\mathcal C}(X,X)_0$ such that $d(1_X)=0$.

The category of $\mathbb{Z}$ graded chain complexes is naturally a dg category.  A dg $\mathcal C$-\emph{module} is a dg enriched functor from $\mathcal C$ to the category of chain complexes. In other words, a $\mathcal C$-module $M$ is the assignment of a chain complex $M(Z)$ to each object $Z$ of $\mathcal C$ together with a $C$-action
\[ \circ: \Hom_{\mathcal C}(Y,Z) \otimes M(Y) \to  M(Z) \] which is associative and unital with respect to the composition in $\mathcal C$.

An important class of dg $\mathcal C$-modules are the \emph{free} of \emph{representable} modules. We say that a $\mathcal C$-module $M$ is \emph{free} or \emph{representable} if there exists a pair $(Y,u)$ which consists of an object $Y$ in $\mathcal C$ and a \emph{universal $0$-cycle} $u \in M(Y)_0$ such that the map
\[ \Hom_{\mathcal C}(Y,Z) \to \ M(Z) \]
induced by $u$ via the module structure is an isomorphism of chain complexes for all $Z$.

\begin{example}
Consider a $k$ algebra $R$ as a dg category with one object $X$, i.e., $\Hom(X,X)=R$, with composition
\[ \circ \colon \Hom (X,X) \otimes_k \Hom (X,X) \to \Hom(X,X), \]
given by multiplicative structure $R \otimes_k R\rightarrow R$. Then the category of dg $R$-modules is the category of chain complexes in $R$.
\end{example}

The motivation for calling a dg category an ``enhancement'' of a triangulated category stems from the following definition.

\begin{definition}\label{pretri}
A dg category $\mathcal C$ is \emph{pretriangulated} if it has a zero object, denoted by $\ast$, such that the following properties hold.
\begin{enumerate}
\item (Closure under shifts)
For an object $X$ in $\mathcal C$ and an $n \in \mathbb{Z}$ the dg $\mathcal C$-module $\Sigma^n \Hom_{\mathcal C}(X,-)$ given by
\[ \left(\Sigma^n \Hom_{\mathcal C}(X,Z) \right)_{n+k} = \Hom_{\mathcal C }(X,Z)_k \] with differential
\[ d(\Sigma^n f) = (-1)^n\cdot \Sigma^n(df) \] is representable.

\item (Closure under cones)
Given a 0-cycle in $\Hom_{\mathcal C}(X,Y)$, the dg $\mathcal C$-module $M$ given by
\[ M(Z)_k = \Hom_{\mathcal C}(Y,Z)_k \oplus \Hom_{\mathcal C}(X,Z)_{k+1}\]
with differential
\[ d(a, b)= (d(a), af-d(b)) \] is representable.
\end{enumerate}
\end{definition}

Underlying any dg category $\mathcal C$ is a preadditive category $\mathcal Z(\mathcal C)$ called the \emph{cycle category}. The category $\mathcal Z(\mathcal C)$ has the same objects as $\mathcal C$, but morphisms are now given by $\Hom_{\mathcal Z(\mathcal C)}(X,Y)=\ker(d \colon \Hom_{\mathcal C}(X,Y)_0\to \Hom_{\mathcal C}(X,Y)_{-1})$, i.e., the morphisms are the $0$-cycles of the complex of morphisms. The \emph{homology category} $H(\mathcal C)$ of a dg category $\mathcal C$ is a quotient of $\mathcal Z(\mathcal C)$. In particular, $H(\mathcal C)$ has the same objects as $\mathcal C$, but morphisms are given by $\Hom_{H(\mathcal C)}(X,Y)=H_0(\Hom_{\mathcal C}(X,Y))$, i.e., morphisms are given by the $0$-th homology groups of the homomorphism complexes. It is the case that if $\mathcal C$ is a pretriangulated dg category, the associated homology category $H(\mathcal C)$ has a natural triangulated structure. A proof of this fact can be found in \cite[\S 3]{bondal-kapranov}, but we describe the shifts and distinguished triangles here for completeness.

Let us assume that $\mathcal C$ is a pretriangulated dg category. A \emph{shift} of an object $X$ in $\mathcal C$ is an object $\Sigma X$ which represents the dg module $\Sigma^{-1}\Hom_\mathcal C(X,-)$ described above. One can take all of the shifts of all objects $X$ of $\mathcal{C}$ and canonically assemble them into an invertible shift functor $X \mapsto \Sigma X$ on $\mathcal C$. This functor induces a shift functor on the homology category $H(\mathcal C)$ (``closure under cones'' in Definition~\ref{pretri}). The distinguished triangles of $H(\mathcal C)$ are triangles that come from mapping cone sequences in $\mathcal C$. More explicitly, a triangle in $\Ho(\mathcal C)$ is \emph{distinguished} if it is isomorphic to the image of a triangle of the form
\[ \xymatrix@1{X \ar[r]^f & Y  \ar[r] & Cf \ar[r] &  \Sigma X}\]
for some $f:X\rightarrow Y$ in $\mathcal C$ (``closure under cones'' in Definition \ref{pretri}).

\begin{example}
Many examples of pretriangulated dg categories come from additive categories, including the pretriangulated hulls of Keller \cite{keller}. In particular, consider the category of modules over a hereditary $k$ algebra $R$. Let $\mathcal{A}= R-\Mod$.  To the additive category $\mathcal A$ we can associate a category of complexes $C(\mathcal A)$ with objects the $\mathbb Z$-graded chain complexes of objects in~$\mathcal A$ and morphisms the chain maps which are homogeneous of degree 0. This category can be made into a dg category $\underline{C}(\mathcal A)$ as follows. Given any two chain complexes $X$ and $Y$ the chain complex of morphisms $\Hom_{\underline{C}(\mathcal A)}(X,Y)$ is given by
\[ \Hom_{\underline{C}(\mathcal A)}(X,Y)_n = \prod_{k \in \mathbb Z} \Hom_\mathcal A(X_k,Y_{k+n}), \]
the abelian group of graded homogeneous morphisms of  degree~$n$. The differential on $\Hom_{\underline{C}(\mathcal A)}(X,Y)$ is given by \[ df = d_Y \circ f - (-1)^n  f \circ d_X \] where $f \in \Hom_{\underline{C}(\mathcal A)}(X,Y)_n$. Composition works as expected.

In this case, the cycle category $\mathcal Z(\underline{C}(\mathcal A))$ is equivalent to the category $C(\mathcal A)$. The homology category $H(\underline{C}(\mathcal A))$ is what is typically called the \emph{homotopy category} $K(\mathcal A)$, which is the category of complexes modulo chain homotopies. We claim that $C(\mathcal A)$ a cofibration category, in the sense to be defined in Section \ref{cofcat}.  Let the class of chain homotopy equivalences be the class of weak equivalences and let the chain maps which are dimension-wise split monomorphisms be the class of cofibrations.
\end{example}

\subsection{The dg orbit category}

There is no reason to assume that a triangulated structure on the orbit category, when it exists, is unique.  However, when it is the triangulated category associated to a dg category, namely, the dg orbit category, it can be regarded as the solution of a universal problem.  Thus, there is a canonical triangulated structure on the orbit category, arising from a dg structure which is unique up to quasi-equivalence.

With this motivation in mind, we give the definition of the dg orbit category.

\begin{definition} \cite{kellerdg}
Let $\mathcal A$ be a dg category and $F \colon \mathcal A \rightarrow \mathcal A$ a dg functor such that $H_0(F)$ is an equivalence.  The \emph{dg orbit category} $\mathcal C$ has the same objects as $\mathcal A$ and morphism complexes defined by
\[ \Hom_{\mathcal C}(X,Y) = \colim_p \bigoplus_{n \geq 0}\Hom_\mathcal A(F^nX, F^pY). \]
\end{definition}

Composition can be defined similarly to the ordinary orbit category, and analogously there is a canonical projection functor $\pi \colon \mathcal A \rightarrow \mathcal C$.  In particular, as categories $H(\mathcal C) \cong H(\mathcal A)/F$.

Again, conditions can be given under which the dg orbit category of a pretriangulated dg category is again a pretriangulated dg category, by showing that it is equivalent to its own dg triangulated hull.

\subsection{Cluster categories} \label{cluster}

The primary example of an orbit category is that of the cluster category, first defined by Buan, Marsh, Reineke, Reiten, and Todorov \cite{bmrrt} as a generalization of a cluster algebra.  Although it can be defined more generally, we consider the specific case of the cluster category associated to an algebra arising from a quiver.

A \emph{quiver} $Q$ is an oriented graph.  We consider here only quivers whose underlying unoriented graph is a Dynkin diagram of type $A$, $D$, or $E$.  (Such graphs have no cycles and are of particular importance in the study of Lie algebras.)  A \emph{representation} of $Q$ over a field $k$ associates to every vertex of $Q$ a $k$-vector space and to every arrow in $Q$ a $k$-linear map. The category of representations of $Q$ over $k$ forms an abelian category $\rep(Q)$.  In homotopy-theoretic language, the \emph{bounded derived category} $\mathcal D^b(Q)$ is the homotopy category of the model category of bounded chain complexes in $\rep(Q)$.  The restrictions we have made on the quiver $Q$ assure that both $\rep(Q)$ and $\mathcal D^b(Q)$ are well-behaved.

\begin{theorem} \cite{happel}
The bounded derived category $\mathcal D^b(Q)$ admits a self-equivalence
\[ \nu \colon \mathcal D^b(Q) \rightarrow \mathcal D^b(Q) \]
such that, for every object $X$, there is an isomorphism of functors
\[ D\Hom(X,-) \rightarrow \Hom(-, \nu X), \]
where $D=\Hom_k(-,k)$.
\end{theorem}

Such a self-equivalence is called a \emph{Serre functor} or \emph{Nakayama functor}.  Additionally, because $\mathcal D^b(Q)$ is a triangulated category, it has an associated shift functor $\Sigma$.

\begin{definition}
The \emph{cluster category} $\mathcal C_Q$ associated to a quiver $Q$ is the orbit category of $\mathcal D^b(Q)$ by the self-equivalence $\nu^{-1}\circ \Sigma^2$.
\end{definition}

In fact, the construction of the cluster category can be placed in to a much more general framework.

\begin{definition}
Let $d$ be an integer. A sufficiently finitary triangulated category $\mathcal T$ is $d$-\emph{Calabi-Yau} if there exists a Serre functor $\nu$ together with a triangulated equivalence $\nu \rightarrow \Sigma^d$.
\end{definition}

From this perspective, we have the following reformulation of the cluster category.

\begin{prop} \cite{keller}
The cluster category $\mathcal C_Q$ is the universal 2-Calabi-Yau category under the bounded derived category $\mathcal D^b(Q)$.
\end{prop}

\section{Topological triangulated categories}

\subsection{Cofibration categories} \label{cofcat}
Topological triangulated categories are defined in terms of cofibration categories. All cofibrantly generated stable model categories satisfy the conditions for a cofibration category.  The following definition is the dual of the one given by Brown for fibration categories \cite[I.1]{brown}; a formulation for cofibration categories can also be found in \cite{schwedecof}.

\begin{definition}
A \emph{cofibration category} is a category $\mathcal C$ equipped with two classes of morphisms, called \emph{cofibrations} and \emph{weak equivalences} which satisfy the following axioms.
  \begin{itemize}
  \item[(C1)] All isomorphisms are cofibrations and weak equivalences.
    Cofibrations are stable under composition.
    The category $\mathcal C$ has an initial object and every morphism from an initial object is a cofibration.
  \item[(C2)] Given two composable morphisms $f$ and $g$ in $\mathcal C$, such that two of the three morphisms $f,g$ and~$gf$ are weak equivalences, then
    so is the third.
  \item[(C3)] Given a cofibration $i \colon A\to B$ and any morphism $f \colon A \rightarrow C$,
    there exists a pushout square
    \begin{equation}
    \begin{aligned}
        \xymatrix{ A \ar[r]^f \ar[d]_i & C \ar[d]^j \\
        B\ar[r] & P}
      \end{aligned}
      \end{equation}
    in $\mathcal C$ and the morphism $j$ is a cofibration.
    If additionally $i$ is a weak equivalence, then so is $j$.
  \item[(C4)] Every morphism in $\mathcal C$ can be factored as
    the composite of a cofibration followed by a weak equivalence.
  \end{itemize}
\end{definition}

We use the term \emph{acyclic cofibration} to denote a morphism that belongs to the class of cofibrations and to the class of weak equivalences. We also note that in a cofibration category a coproduct $B\vee C$ of any two objects in $\mathcal C$ exists.  The canonical morphisms from $B$ and $C$ to $B\vee C$ are cofibrations. The \emph{homotopy category} of a cofibration category is a localization at the class of weak equivalences, i.e., a functor $\gamma \colon \mathcal C \to \Ho(\mathcal C)$ that takes all weak equivalences to isomorphisms which is initial among such functors.

If one prefers to work in model categories, one can obtain a cofibration category by restricting to the full subcategory of cofibrant objects and forgetting the fibrations. For the purposes of this article we primarily consider examples which arise from model categories.

Like in the case with dg categories, cofibration categories which satisfy some extra conditions are enhancements of triangulated categories. A cofibration category is \emph{pointed} if every initial object is also terminal. We denote this \emph{zero object} by $\ast$. In a pointed cofibration category, the axiom $(C4)$ provides a \emph{cone} for every object $A$, i.e., a cofibration $i_A \colon A\to CA$ whose target is weakly equivalent to $\ast$.

Given a pointed cofibration category $\mathcal C$, the \emph{suspension} $\Sigma A$ of an object $A$ in $\mathcal C$ is the quotient of the \emph{cone inclusion}. This is equivalent to a pushout
\[ \xymatrix{ A \ar[r]^-{i_A} \ar[d]& CA \ar[d] \\
\ast \ar[r] & \Sigma A. } \]
As with pretriangulated dg categories, one can assemble the suspension construction into a functor $\Sigma \colon \Ho(\mathcal C)\to \Ho(\mathcal C)$ on the level of homotopy categories.  For cofibrations, an argument is given in the Appendix to \cite{schwedecof}.

The class of cofibrations in $\mathcal C$ allow us to define distinguished triangles in $\Ho(\mathcal C)$. In particular, each cofibration $j \colon A \rightarrow B$ in a pointed cofibration category $\mathcal C$ gives rise to a natural \emph{connecting morphism} $\delta(j) \colon B/A \to \Sigma A$ in $\Ho(\mathcal C)$. The \emph{elementary distinguished triangle} induced by the cofibration $j$ is the triangle
\[ \xymatrix@1{A \ar[r]^j & B \ar[r]^q &  B/A \ar[r]^{\delta(j)} & \Sigma A} \]
where $q \colon B\to B/A$ is a \emph{quotient morphism}. A \emph{distinguished triangle} is any triangle that is isomorphic to the elementary distinguished triangle of a cofibration in the homotopy category.

A pointed cofibration category is \emph{stable} if the suspension functor $\Sigma \colon \Ho(\mathcal C) \to \Ho(\mathcal C)$ is a self-equivalence. The suspension functor and the class of distinguished triangles make the homotopy category $\Ho(\mathcal C)$ into a triangulated category.

\begin{definition}
A triangulated category is \emph{topological} if it is equivalent, as a triangulated category, to the homotopy category of a stable cofibration category.
\end{definition}

The adjective ``topological" does not imply that the category or its hom-sets have a topology, but rather that these examples are constructed by methods in the spirit of abstract homotopy theory.

\subsection{Topological triangulated categories arising from algebraic ones}

One can demonstrate that the cycle category $\mathcal Z(\mathcal B)$ of a pretriangulated dg category $\mathcal B$ is a cofibration category. A closed morphism is a \emph{weak equivalence} if it becomes an isomorphism in the homology category. A closed morphism $i \colon A \rightarrow B$ is a \emph{cofibration} if:
\begin{itemize}
\item the induced chain morphism $\Hom_\mathcal B(i,Z)$ is surjective for every object $Z$ of $\mathcal B$ and
\item the kernel $\mathcal B$-module
\[ Z \mapsto \ker \left[ \Hom_\mathcal B(i,Z) \colon \Hom_\mathcal B(B,Z) \rightarrow \Hom_\mathcal B(A,Z) \right] \]
is representable.
\end{itemize}

Notice that given that the module $(C,u)$ represents the kernel of $\mathcal B(i,-)$, then, by definition, there exists a universal $0$-cycle $u \colon B \rightarrow C$ such that for every $Z$ of $\mathcal B$  the sequence of cycle groups
\[ \xymatrix@1{0 \ar[r] & \Hom_{\mathcal Z(\mathcal B)}(C,Z) \ar[r]^{u^*} & \Hom_{\mathcal Z(\mathcal B)}(B,Z) \ar[r]^{i^*} & \Hom_{\mathcal Z(\mathcal B)}(A,Z)} \]
is exact. In particular, $u \colon B \rightarrow C$ is a cokernel of $i \colon A \rightarrow B$ in the category $\mathcal Z(\mathcal B)$.

The following proposition is due to Schwede.

\begin{prop}\label{schwede1} \cite[3.2]{schwedepre}
Let $\mathcal B$ be a pretriangulated dg category.  Then the cofibrations and weak equivalences make the cycle category $\mathcal Z(\mathcal B)$ into a stable cofibration category in which every object is fibrant. Moreover, the homotopy category $\Ho(\mathcal Z(\mathcal B))$ is equivalent, as a triangulated category, to the homology category $H_0(\mathcal B)$. In particular, every algebraic triangulated category is a topological triangulated category.
\end{prop}

\subsection{Homotopy colimits of cofibration categories}

We would like to define topological orbit categories via a coequalizer construction.  Therefore, we require a notion of homotopy colimits of cofibration categories.

\begin{definition}
Let $\mathcal D$ be a small category, and $\mathcal M$ a $\mathcal D$-shaped diagram of functors between cofibration categories $F_{\alpha, \beta}^\theta \colon \mathcal M_\alpha \rightarrow \mathcal M_\beta$.  (Here the superscript $\theta$ allows us to distinguish between different arrows $\alpha \rightarrow \beta$ in $\mathcal D$.) Then the \emph{homotopy colimit} of $\mathcal M$, denoted by $\mathcal Colim_\alpha\mathcal M_\alpha$, is defined to be the category obtained from the disjoint union of the model categories in $\mathcal M$ by inserting weak equivalences $x_\beta \rightarrow x_\alpha$ between objects $x_\alpha$ in $\mathcal M_\alpha$ and $x_\beta$ in $\mathcal M_\beta$ if there exists a weak equivalence $F_{\alpha, \beta}^\theta(x_\alpha) \rightarrow x_\beta$ in $\mathcal M_\beta$.  We further assume that, if such a weak equivalence already exists (in the case where $\alpha=\beta$), we do not add an additional one, and that we impose the appropriate relation on composites: if there exist two weak equivalences $F_{\alpha, \beta}^\theta(x_\alpha) \rightarrow x_\beta$ and $F_{\beta, \gamma}^\psi(x_\beta) \rightarrow x_\gamma$, then the two possible ways of obtaining weak equivalences $x_\alpha \rightarrow x_\gamma$ are identified.
\end{definition}

The definition of a homotopy colimit of (stable) model categories is given in \cite{hocolim}, but the definition can be modified as above more general (stable) cofibration categories.  In fact, the homotopy colimit of a diagram of model categories is not generally still a model category, but a more general homotopy theory.  One can take it to be a cofibration category. The question, then, as to whether an orbit category is triangulated can now be understood as a question of whether or not a homotopy coequalizer of cofibration categories is a stable cofibration category.

\subsection{Topological orbit categories}

We now define a topological orbit category as a generalization of a dg orbit category. As a consequence of Proposition \ref{schwede1}, this definition includes all of the previously known algebraic examples.  Let $\mathcal T$ be a stable cofibration category and $F \colon T \rightarrow T$ a standard equivalence of cofibration categories, i.e., a functor inducing a triangulated equivalence on the homotopy category.  In the case where $\mathcal T$ is a stable model category, then we ask that $F$ be one of the adjoint maps in a derived Morita equivalence (see \cite{schwedeshipley}).

\begin{definition}
The \emph{topological orbit category} $\mathcal T/F$ is the homotopy coequalizer of the diagram
\[ \xymatrix@1{\mathcal T \ar@<.5ex>[r]^\id \ar@<-.5ex>[r]_F & \mathcal T \ar[r] & \mathcal T/F.} \]
\end{definition}

\subsection{Topological cluster categories}

We can now consider our motivating example of cluster categories from the point of view of triangulated orbit categories.  For simplicity, we work in the context of stable model categories.

\begin{definition}
A \emph{Serre functor} on a stable model category $\mathcal T$ is a Quillen functor $\nu \colon \mathcal T \rightarrow \mathcal T$ inducing a Serre functor of triangulated categories on $\Ho(\mathcal T)$.
\end{definition}

\begin{definition}
Let $\mathcal T$ be a stable model category which admits a Serre functor.  The \emph{cluster category} of $\mathcal T$ is the topological orbit category of $\mathcal T$ by the self-equivalence $\nu \circ \Sigma^2$.
\end{definition}

While orbit categories of a purely topological origin are to appear in \cite{long}, we include the following motivating example.

\begin{example}
Let $Q$ be the quiver $A_{m}$.  For a field $k$, let $A$ be the path algebra $kQ$.  This example is well-studied from the perspective of cluster categories, and here we extend to our topological perspective using the framework of ring spectra and stable model categories.  We would like to show that the stable model category approach recovers the original construction, under base change along Quillen equivalences.  Hence, topological cluster categories recover this known construction.

Work of Shipley shows that we may consider $A$ as a ring spectrum under the image of the Eilenberg-Mac Lane functor $H$ \cite{ship}. Note that if $k$ is a commutative ring and $A$ is a $k$-algebra, then $HA$ is an $Hk$-algebra and $Hk$ is a commutative ring spectrum.  The category $\Mod_{HA}$ of $HA$-modules, which we denote by $\mathcal T$ for simplicity, forms a stable model category by \cite[3.1.1]{schwedeshipley}.  As such, $\Ho(\mathcal{T})$ is a triangulated category.

Recall from Section \ref{cluster} that the cluster category is the triangulated orbit category of $\mathcal D^{b}(A)$ under the action of the triangulated equivalence $F$ given by $M \mapsto \tau \Sigma M$ where $\tau$ is the Auslander-Reiten translation of $\mathcal D^{b}(A)$ and $\Sigma$ is the suspension functor.    We want to construct the orbit category of $\Ho(\mathcal T)$ by $F$.  In this case, however, $F^{-1}$ is more easily described; it is given by $M \mapsto \Sigma^{2}\nu$, where $\nu$ is the Nakayama functor
\[ \nu = - \otimes_{A}^{\mathbb{L}} \Hom_{k}(A,k). \]

We can now explicitly construct the orbit category $\Ho(\mathcal{T})/F^{-1}$, showing that it is given by the triangulated category $\Ho(\Mod_{HB})$ for $B$ a differential graded algebra and $HB$ its corresponding ring spectrum under the Eilenberg-Mac Lane functor. Consider the differential graded algebra
\[ B=A\oplus \Sigma^{3}\Hom_{k}(A,k) \]
where the $k$-algebra $A$ is in degree $0$ and the bimodule $\Sigma^{3}\Hom_{k}(A,k)$ is in degree $3$; $B$ has trivial differential. Keller proves that $\mathcal D^{b}(B)$ is the orbit category $\mathcal D^{b}(A)/F^{-1}$ \cite{keller}. We want to show that this equivalence holds on the level of stable model categories, which requires a somewhat subtle equivariance problem.

For his construction, Keller uses some facts about group actions on triangulated categories, specifically that $\mathcal D^{b}(B)$ admits a canonical action by the braid group on $n+1$ strings. This action was investigated by Khovanov and Seidel in \cite{ks} and independently in a similar context by Rouquier and Zimmermann \cite{rz}.  Khovanov and Seidel write down explicit triangulated self-equivalences $R_{1},\ldots,R_{m}$ of $\mathcal D^{b}(B)$ which satisfy relations
\begin{itemize}
    \item $R_{i}R_{i+1}R_{i}\cong R_{i+1}R_{i}R_{i+1}$ for $1\leq i< m$;
    \item $R_{j}R_{k}\cong R_{k}R_{j}$ for $|j-k|\geq 2$
\end{itemize}
which are precisely the defining relations of the braid group $B_m$.  In other words, these triangulated self equivalences $R_{i}$ give an action of $B_{m+1}$ on $\mathcal D^{b}(B).$

Each of canonical generators of the braid group action $R_{i}$ is given by a self-equivalence $R_{i}'\otimes_{A} \Sigma^{2}\Hom_{k}(A,k)$. Moreover, each of these $R_{i}'$ is given by a complex of bimodules which has $A$ in degree $0$ and some product of projectives in degree $1$.  There is a morphism of triangulated functors $\varphi_{i} \colon R_{i}'\rightarrow 1$, where $1$ is the self equivalence of $\mathcal D^{b}(B)$ which shifts degree of the dga by one.

The cone on each $\varphi_{i}$ belongs to $\per(B)$, which is equal to the smallest triangulated subcategory of itself stable under direct factors and containing these cones. Thus, the action of $B_{m+1}$ becomes trivial in $\mathcal D^{b}(B)/\per(B)$ and in a certain sense, this quotient is the largest one in which the $\varphi_{i}$ become invertible.

If our $\Ho(\mathcal{T})/F^{-1}$ is to be equivalent, via base change along Quillen equivalences, to $\mathcal D^{b}(B)$, we need to show how this canonical braid group action translates under Quillen equivalence. In this particular example, the canonical action is given by triangulated equivalences that come from tensoring with bimodules, so the extension is straightforward.  One can check that the triangulated self equivalences $R_{i}$ of Khovanov and Seidel are still triangulated equivalences $R_{i}'\otimes_{A} \Sigma^{2}\Hom_{Hk}(HA,Hk)$, where now $\Hom_{Hk}(HA,Hk)$ is the bimodule given by the mapping spectrum as in Dwyer-Greenlees-Iyengar \cite{dgi}.  In this case, each of the $R_{i}'$ is given by a bimodule over a symmetric ring spectrum which, using the terminology of Dwyer-Greenlees-Iyengar, is built from $HA$.  Therefore, we still obtain maps $\varphi_{i} \colon R_{i}'\rightarrow 1$ whose cones are perfect. A calculation shows that $\per(HB)$ is equal to its smallest triangulated subcategory stable under direct factors and containing the necessary cones, and that the action of $B_{m+1}$ becomes trivial on the orbit category $\Ho(\mathcal{T})/F^{-1}$.  Base change along Quillen equivalences now recovers the original algebraic result.
\end{example}

\end{document}